\documentclass[11pt]{article}
\usepackage{epic,latexsym,amssymb}
\usepackage{color}
\usepackage{tikz}
\usepackage{amsfonts,epsf,amsmath}

\newtheorem{thm}{Theorem}

\newtheorem{cor}{Corollary}

\newcommand{\qed}{$\Box$}

\newcommand{\smallqed}{{\tiny ($\Box$)}}

\newcommand{\cF}{{\cal F}}

\newcommand{\proof}{\noindent\textbf{Proof. }}

\let\oldenumerate\enumerate
\renewcommand{\enumerate}{
  \oldenumerate
  \setlength{\itemsep}{0pt}
  \setlength{\parskip}{0pt}
  \setlength{\parsep}{0pt}
}

\begin{document}

\title{Inequalities between Partial Domination and Independent Partial Domination in Graphs}

\author{Odile Favaron and Pawaton Kaemawichanurat$^{1}$
\\ \\
$^{1}$Theoretical and Computational Science Center \\
and Department of Mathematics, Faculty of Science,\\
King Mongkut's University of Technology Thonburi, \\
Bangkok, Thailand\\
\small \tt Email: odile.favaron@lri.fr; pawaton.kae@kmutt.ac.th}

\date{}
\maketitle

\begin{abstract}
For a graph $G$, a vertex subset $S \subseteq V(G)$ is said to be $K_{k}$-isolating if $G - N_{G}[S]$ does not contain $K_{k}$ as a subgraph. The $K_{k}$-isolation number of $G$, denoted by $\iota_{k}(G)$, is the minimum cardinality of a $K_{k}$-isolating set of $G$. Analogously, $S$ is said to be independent $K_{k}$-isolating if $S$ is a $K_{k}$-isolating set of $G$ and $G[S]$ has no edge. The independent $K_{k}$-isolation number of $G$, denoted by $\iota'_{k}(G)$, is the minimum cardinality of an independent $K_{k}$-isolating set of $G$. A vertex subset $D \subseteq V(G)$ is said to be dominating if $V(G) \setminus N_{G}[S] = \emptyset$. Moreover, if $G[D]$ has no edge, then $D$ is an independent dominating set. The cardinality of a smallest dominating set is the domination number and is denoted by $\gamma(G)$, similarly, the cardinality of a smallest independent dominating set is the independent domination number and is denoted by $i(G)$. Clearly, when $k = 1$, we have $\gamma(G) = \iota_{1}(G)$ and $i(G) = \iota'_{1}(G)$. For classic results between $\gamma(G)$ and $i(G)$, in 1978, Allan and Laskar proved that $\gamma(G) = i(G)$ for all $K_{1, 3}$-free graphs and this result was generalized to $K_{1, r}$-free graphs by Bollob$\acute{a}$s and Cockayne in 1979. In 2013, Rad and Volkmann proved that the ratio $i(G)/\gamma(G)$ is at most $\Delta(G)/2$ when $\Delta(G) \in \{3, 4, 5\}$. Further, Furuya et. al. proved that when $\Delta(G) \geq 6$, we have $i(G)/\gamma(G) \leq \Delta(G) - 2\sqrt{\Delta(G)} + 2$. In this paper, for a smallest $K_{k}$-isolating set $S$, we prove that $\iota'_k(G)\le -\frac{\iota_k^2(G)}{\ell} +i_k(G)(\Delta +2)-\ell \Delta$
where $\ell$ is the number of some specific vertices of $S$ such that the union of their closed neighborhoods in $S$ is $S$. We prove that this bound is sharp. A special case of our main theorem implies $\iota'_{k}(G)/\iota_{k}(G) \leq \Delta(G) - 2\sqrt{\Delta(G)} + 2$. Further, we find an inequality between $\iota'_{k}(G)$ and $\iota_{k}(G)$ when $G$ is $K_{1, r}$-free graph. This also generalizes the result of Bollob$\acute{a}$s and Cockayne.
\end{abstract}

{\small \textbf{Keywords:} Partial-domination.} \\
\indent {\small \textbf{AMS subject classification:} 05C69}

\section{Introduction and Background}
Let $G$ be a simple graph with vertex set $V(G)$ and edge set $E(G)$ of order~$n(G) = |V(G)|$ and size $m(G) = |E(G)|$. We denote the \emph{degree} of $v$ in $G$ by $deg_G(v)$ and denote the maximum degree of $G$ by $\Delta(G)$. A \emph{neighbor} of a vertex $v$ in $G$ is a vertex $u$ which is adjacent to $v$. The \emph{open neighborhood} $N_G(v)$ of a vertex $v$ in $G$ is the set of neighbors of $v$. That is $N_G(v) = \{u \in V(G) \, | \, uv \in E(G)\}$. The \emph{closed neighborhood} of $v$ is $N_G[v] = N_G(v) \cup \{v\}$. For a subset $S \subseteq V(G)$, we use $N_{S}(v)$ to denote $N_{G}(v) \cap S$ and $deg_{S}(v) = |N_{G}(v) \cap S|$, moreover, we use $N_{S}[v]$ to denote $N_{G}[v] \cap S$. The \emph{neighborhood} of a vertex subset $S$ of $G$ is the set $N_G(S) = \cup_{v \in S} N_G(v)$. The \emph{closed neighborhood} of $S$ in $G$ is the set $N_G[S] = N_G(S) \cup S$. The subgraph of $G$ induced by $S$ is denoted by $G[S]$. The subgraph obtained from $G$ by deleting all vertices in $S$ and all edges incident with vertices in $S$ is denoted by $G - S$. The \emph{distance} between two vertices $u$ and $v$ in a connected graph $G$ is the length of a shortest $(u,v)$-path in $G$ and is denoted by $d_G(u,v)$. We denote the \emph{clique} on $n$ vertices by $K_n$. A \emph{star} $K_{1, n}$ is a graph of $n + 1$ vertices obtained by joining $n$ vertices to one vertex. A graph $G$ is $H$\emph{-free} if $G$ does not contain $H$ as an induced subgraph.
\vskip 5 pt

A vertex subset $S$ of a graph $G$ is a \emph{dominating set} of $G$ if every vertex in $V(G) \setminus S$ is adjacent to a vertex in $S$. The cardinality of a smallest dominating set of $G$ is called the \emph{domination number} of $G$ and is denoted by $\gamma(G)$. Moreover, $S$ is an \emph{independent dominating set} of $G$ if $S$ is a dominating set of $G$ and there is no edge in $G[S]$. The cardinality of a smallest independent dominating set of $G$ is called the \emph{independent domination number} of $G$ and is denoted by $i(G)$.
\vskip 5 pt

Recently, Caro and Hansberg~\cite{CaHa17} generalized the concept of domination by focusing on a vertex subset $S \subseteq V(G)$ so that $G - N_{G}[S]$ contains no forbidden subgraph. Let $G$ be a graph and $\cF$ a family of graphs. A vertex subset $S \subseteq V(G)$ is said to be $\cF$\emph{-isolating} if $G - N_{G}[S]$ does not contain $H$ as a subgraph for all $H \in \cF$. Obviously, when $\cF = \{K_{1}\}$, $S$ is a dominating set. The $\cF$-\emph{isolation number} of $G$, denoted by $\iota_{\cF}(G)$, is the minimum cardinality of an $\cF$-isolating set of $G$.  Analogously, $S$ is said to be \emph{independent} $\cF$\emph{-isolating} if $S$ is an $\cF$-isolating set of $G$ and $G[S]$ has no edge. The \emph{independent} $\cF$-\emph{isolation number} of $G$, denoted by $\iota'_{\cF}(G)$, is the minimum cardinality of an independent $\cF$-isolating set of $G$. 
In this paper, we consider a family $\cF$ reduced to a clique $K_k$ for a positive integer $k$ and  we denote $\iota_{\{K_{k}\}}(G)$ and $\iota'_{\{K_{k}\}}(G)$ by $\iota_{k}(G)$ and $\iota'_{k}(G)$. A smallest $\{K_{k}\}$-isolating set is called an $\iota_{k}$-\emph{set} and
a smallest independent $\{K_{k}\}$-isolating set is called an $\iota'_{k}$-\emph{set}. We generalize to $\iota_k(G)$ and $\iota'_k(G)$ two known results related to $\iota_1(G)=\gamma(G)$ and $\iota'_1(G)=i(G)$.
\vskip 5 pt

For classic results between the domination number and the independent domination number, Allan and Laskar~\cite{AL 1978} proved that $\gamma(G) = i(G)$ for all $K_{1, 3}$-free graphs and this result was generalized to $K_{1, r}$-free graphs by Bollob$\acute{a}$s and Cockayne~\cite{BoCo1979}. That is :

\begin{thm}\label{thm bo}\cite{BoCo1979}
Let $G$ be a $K_{1, r}$-free graph where $r \geq 3$. Then $i(G) \leq (r - 2)(\gamma(G) - 1) + 1$.
\end{thm}

\indent In 2013, Rad and Volkmann proved that the ratio $i(G)/\gamma(G) \leq \Delta(G)/2$ when $3 \leq \Delta(G) \leq 5$. When $\Delta(G) \geq 6$, they conjectured analogously that $i(G)/\gamma(G) \leq \Delta(G)/2$. However, the conjecture was disproved by Furuya et. al. \cite{FOS} with the upper bound $\Delta(G) - 2\sqrt{\Delta(G)} + 2$ sharp for every $\Delta$ equal to a square. That is :

\begin{thm}\label{thm fu}\cite{FOS}
For a graph $G$, $i(G)/\gamma(G) \leq \Delta(G) - 2\sqrt{\Delta(G)} + 2$.
\end{thm}
\vskip 5 pt

\section{Main results}
In this section, we state our main results and prove that all these results are sharp. It is worth noting that an $\iota'_{k}$-set is a $\{K_k\}$-isolating set. By the minimality, $\iota_{k}(G) \leq \iota'_{k}(G)$ for any graph $G$. In the following, when no ambiguity can occur, we write $\iota_{k}, \iota'_{k}$ and $\Delta$ rather than $\iota_{k}(G), \iota'_{k}(G)$ and $\Delta(G)$, respectively.
We prove that :
\vskip 5 pt

\begin{thm}\label{thm 1}
Let $G$ be a graph with maximum degree $\Delta$ and let $S$ be an $\iota_k(G)$-set for some positive integer $k$. Let $v_1,v_2, ,
\cdots  v_{\ell}$ be a sequence of vertices of $S$ such that $v_1$ has minimum degree in $S$ and recursively $v_{i+1}$ has minimum degree in $S\setminus N_S[\{v_1,v_2,\cdots ,v_i\}]$ until $S=\cup _{i=1}^{\ell}N_S[v_i]$. Then

\begin{center}
$\iota'_k(G)\le -\frac{\iota_k^2(G)}{\ell} +i_k(G)(\Delta +2)-\ell \Delta.$
\end{center}
\end{thm}
\vskip 5 pt

\noindent We will give the proof in Section \ref{sec 3}.
The following corollary of Theorem \ref{thm 1}   generalizes Theorem \ref{thm fu} to all positive values of $k$.
\vskip 5 pt

\begin{cor}\label{cor 1}
For a graph $G$ and an integer $k \geq 1$,  $\iota'_{k}(G)/\iota_{k}(G) \leq \Delta - 2\sqrt{\Delta} + 2$.
\end{cor}
\proof
The maximum of the function
 $f : \mathcal{R}^{+} \rightarrow \mathcal{R}$ defined by
\begin{center}
$f(x)= -\frac{\iota_k^2}{x} +i_k(\Delta +2)-x\Delta$
\end{center}
is attained when $x=\frac{\iota _k}{\sqrt{\Delta}}$ and is equal to
$-\sqrt{\Delta}\iota _k+(\Delta +2)\iota _k-\sqrt{\Delta}\iota _k$.
By Theorem \ref{thm 1},
\begin{center}
$\iota '_k\le f(\frac{\iota_k}{\sqrt{\Delta}})=(\Delta - 2\sqrt{\Delta} + 2)\iota_k$.
\end{center}
\begin{center}
Hence $~~\iota'_{k}(G)/\iota_{k}(G) \leq \Delta - 2\sqrt{\Delta} + 2$.
\end{center}

\qed
\vskip 5 pt

\indent In our last main result, we generalize Theorem \ref{thm bo} by establishing the upper bound of $\iota'_{k}$ in terms of $\iota_{k}$ and $r$ in $K_{1, r}$-free graphs. We find the same upper bound as that of Theorem \ref{thm bo}. The proof is provided in Section \ref{sec 4}.
\vskip 5 pt

\begin{thm}\label{thm 2}
Let $G$ be a $K_{1, r}$-free graph where $r \geq 3$. Then $\iota'_{k} \leq (r - 2)(\iota_{k} - 1) + 1$.
\end{thm}
\vskip 5 pt

\indent We conclude this section by giving a construction of graphs satisfying the equality in Theorems \ref{thm 1}, \ref{thm 2} and Corollary \ref{cor 1}.
\vskip 5 pt

\noindent \textbf{The graphs $G(t, s)$}\\
\indent Let $s, t, k$ be positive integers such that $s + t - 1 \geq k$. For $1 \leq i \leq t$, we let $K^{i, 1}_{k}, K^{i, 2}_{k}, ..., K^{i, s}_{k}$ be $s$ disjoint copies of a clique $K_{k}$. Let $x_{1}, x_{2}, ..., x_{t}$ be $t$ vertices. The graph $G(t, s)$ is obtained from $K^{i, 1}_{k}, K^{i, 2}_{k}, ..., K^{i, s}_{k}$ for $1 \leq i \leq t$ and $x_{1}, x_{2}, ..., x_{t}$ by joining each $x_{i}$ to a vertex of $K^{i, i'}_{k}$, $y_{i, i'}$ say, for all $1 \leq i' \leq s$ and form $x_{1}, x_{2}, ..., x_{t}$ a clique. Observe that $deg_{G(t, s)}(x_{i}) = \Delta(G(t, s)) = s + t - 1 \geq k$.
\vskip 2 pt

\indent We see that $\{x_{1}, x_{2}, ..., x_{t}\}$ is a $\{K_{k}\}$-isolating set of $G(t, s)$. Thus, $\iota_{k}(G(t, s)) \leq t$. Let $S$ be an $\iota_{k}$-set of $G(t, s)$. To be adjacent to cliques $K^{i, 1}_{k}$, we have that $(\{x_{i}\} \cup V(K^{i, 1}_{k})) \cap S \neq \emptyset$. Thus, $\iota_{k}(G(t, s)) = |S| \geq t$ implying that $\iota_{k}(G(t, s)) = t$.
\vskip 2 pt

\indent Now, we let $\iota'_{k}(G(t, s)) = t'$. We will show that $t' = s(t - 1) + 1$.  Clearly, $\{x_{1}\} \cup \{y_{i, i'} : 2 \leq i \leq t$ and $1 \leq i' \leq s\}$ is an independent $\{K_{k}\}$-isolating set of $G(t, s)$. So, $\iota'_{k}(G(t, s)) \leq s(t - 1) + 1$. Let $S'$ be an $\iota'_{k}$-set of $G(t, s)$. By the independence of $S'$, $|S' \cap \{x_{1}, x_{2}, ..., x_{t}\}| \leq 1$. Without loss of generality, we let $x_{1}, x_{2}, ..., x_{t - 1} \notin S'$. Hence, $S' \cap V(K^{i, i'}_{k}) \neq \emptyset$ for all $1 \leq i \leq t - 1$ and $1 \leq i' \leq s$. Moreover, to be adjacent to $K^{t, 1}_{k}$, we have that $S' \cap (V(K^{t, 1}_{k}) \cup \{x_{t}\}) \neq \emptyset$. Hence, $\iota'_{k}(G(t, s)) = |S'| \geq s(t - 1) + 1$ implying that $t' = s(t - 1) + 1$.
\vskip 5 pt

\indent If we let $s=t^2-t+1$ with $t^2\ge k$,  then $\Delta(G(t, s)) = t^{2}$ and
$t' = s(t - 1) + 1 = t^3 -2t^2 +2t .$

Hence
$t'/t = \Delta - 2\sqrt{\Delta} + 2$
and
$t'= -t^2 +t(\Delta +2)-\Delta.$

This shows that the bounds of Corollary \ref{cor 1} and of Theorem \ref{thm 1} in the case $\ell =1$ are attained by arbitrarily large graphs.
\vskip 10 pt

\indent We can construct a graph $\tilde{G}$ satisfying the equality for the bound in Theorem \ref{thm 1} for any positive value of $\ell$ by letting $\tilde{G}$ be the disjoint union of $G_{1}, ..., G_{\ell}$ where each $G_{i}$ is a copy of $G(t, t^{2} - t + 1)$ as defined in the above paragraph. Similarly, we have $\Delta(\tilde{G}) = t^2, \iota_{k}(\tilde{G}) = t\ell$ and $\iota'_{k}(\tilde{G}) = (t^{3} - 2t^{2} + 2t)\ell$. Hence,
\begin{align}
\iota'_{k}(\tilde{G}) &= (t^{3} - 2t^{2} + 2t)\ell \notag\\
                  & §= -\frac{t^2\ell ^2}{\ell} +t\ell(t^2+2)-\ell t^2\notag\\
                     &= -\frac{(\iota _k(\tilde{G}))^2}{\ell} + \iota_k(\tilde{G})(\Delta(\tilde{G})+2) - \ell \Delta (\tilde{G}). \notag
%
\end{align}
The graph $\tilde{G}$ is not connected. We can make it connected when $k \geq 3$ by joining with a path of length at least four one vertex of $K^{t, s}_{k}$ of $G_{j}$ to one vertex of $K^{1, 1}_{k}$ of $G_{j + 1}$ for $1\le j\le \ell -1$.
For the resulting graph $\hat{G}$, $\Delta(\hat{G}) = t^2, \iota_{k}(\hat{G}) = t\ell$ and $\iota'_{k}(\hat{G}) = (t^{3} - 2t^{2} + 2t)\ell$. Hence, $\iota'_{k}(\hat{G})$ satisfies the bound in Theorem \ref{thm 1}.
\vskip 5 pt

\indent Finally, if we let $s = r - 2$, the graph $G(t, s)$ is $K_{1,r}$-free and $t' = s(t - 1) + 1 = (r - 2)(t - 1) + 1.$ This shows that the bound of Theorem \ref{thm 2} is sharp.
\vskip 10 pt

\section{Proof of Theorem \ref{thm 1}}\label{sec 3}
First, we restate Theorem \ref{thm 1}.
\vskip 2 pt

\noindent \textbf{Theorem \ref{thm 1}}
Let $G$ be a graph with maximum degree $\Delta$ and let $S$ be an $\iota_k(G)$-set for some positive integer $k$. Let $v_1,v_2, ,
\cdots  v_{\ell}$ be a sequence of vertices of $S$ such that $v_1$ has minimum degree in $S$ and recursively $v_{i+1}$ has minimum degree in $S\setminus N_S[\{v_1,v_2,\cdots ,v_i\}]$ until $S=\cup _{i=1}^{\ell}N_S[v_i]$. Then

\begin{center}
$\iota'_k(G)\le -\frac{\iota_k^2(G)}{\ell} +i_k(\Delta +2)-\ell \Delta.$
\end{center}


\proof
Let $\{v_1,v_2,\cdots,v_{\ell}\}$ be a sequence of vertices of $S$ as defined in the theorem.
Initially, we let $S_{0} = S$. Then, we let for $1\le i\le \ell$,



\begin{center}
$S_{i} = S_{i-1} \backslash N_{S_{i-1}}[v_{i}]$ and $N_{S_{i-1}}(v_{i}) = \{v_{i}^{1}, v_{i}^{2}, ...,  v_{i}^{j_{i}}\}$
\end{center}
\vskip 2 pt

\noindent where $j_{i} = deg_{S_{i-1}}(v_{i})$. It is worth noting that


\begin{align}
\emptyset = S_{\ell} \subset S_{\ell - 1} \subset  S_{\ell - 2} \subset \cdot \cdot \cdot \subset S_{1} \subset S_{0}.\notag
\end{align}

\noindent \textbf{Claim 1 :} $\{v_{1}, v_{2}, ..., v_{\ell}\}$ is an independent set.\\
\proof This is a consequence of the construction of the sequence $v_1,\cdots ,v_{\ell}$ since for $2\le j\le \ell$, $v_j\notin \bigcup _{i=1}^{j-1}N_{S_{i-1}}[v_i]$. \smallqed

\vskip 5 pt

\noindent \textbf{Claim 2 :} $\Sigma^{\ell}_{i = 1}(deg_{S_{i - 1}}(v_{i}) + 1) = |S|$, in particular, $\Sigma^{\ell}_{i = 1} {deg}_{S_{i - 1}}(v_{i}) = |S| - \ell$.\\
\proof Clearly $\bigcup^{\ell}_{i = 1} N_{S_{i - 1}}[v_{i}] = S$ and $N_{S_{i - 1}}[v_{i}] \cap N_{S_{j - 1}}[v_{j}] = \emptyset$. Thus $\Sigma^{\ell}_{i = 1}|N_{S_{i - 1}}[v_{i}]| = |S|$. Because $|N_{S_{i - 1}}[v_{i}]| = {deg}_{S_{i - 1}}(v_{i}) + 1$, it follows that $\Sigma^{\ell}_{i = 1}(deg_{S_{i - 1}}(v_{i}) + 1) = |S|$. Hence, $\Sigma^{\ell}_{i = 1} deg_{S_{i - 1}}(v_{i}) = |S| - \ell$. This completes the proof.
\smallqed
\vskip 5 pt

\indent For a clique $K_{k}$ and a vertex $v \in V(G)$, we say that $K_{k}$ \emph{is adjacent to} $v$ (or vice versa) if $v$ is adjacent to a vertex of $K_{k}$ or is a vertex of $K_k$.
\vskip 5 pt
Let $A=N_{G\setminus S}(S) - N_{G\setminus S}(v_1,v_2,\cdots , v_{\ell})$.

\noindent \textbf{Claim 3 :} $|A|\le \Sigma^{\ell}_{i = 1}deg_{S_{i - 1}}(v_{i})(\Delta - deg_{S_{i - 1}}(v_{i})).$


\proof Clearly ${deg}_{S_{i - 1}}(v^{j}_{i}) + {deg}_{G \backslash S_{i - 1}}(v^{j}_{i}) = {deg}_{G}(v^{j}_{i}) \leq \Delta$. Thus, from the choice of $v_i$,
\begin{align}
{deg}_{G \backslash S_{i - 1}}(v^{j}_{i}) \leq \Delta - {deg}_{S_{i - 1}}(v^{j}_{i}) \leq \Delta - deg_{S_{i - 1}}(v_{i})\notag
\end{align}
\noindent for $1\le j\le j_i$. Therefore $v^{j}_{i}$ has at most $ \Delta - {deg}_{S_{i - 1}}(v_{i})$ neighbors in $G\setminus S$.
Hence, and since $j_{i} = deg_{S_{i - 1}}(v_{i})$,
\begin{align}
|A| &\le \Sigma^{\ell}_{i = 1} \Sigma^{j_{i}}_{j = 1} {deg}_{G \setminus S}(v^j_{i})\notag\\
    &\le \Sigma^{\ell}_{i = 1} \Sigma^{j_{i}}_{j = 1} (\Delta - {deg}_{S_{i - 1}}(v_{i}))\notag\\
    &= \Sigma^{\ell}_{i = 1}deg_{S_{i - 1}}(v_{i})(\Delta - deg_{S_{i - 1}}(v_{i}))\notag.
\end{align}
\smallqed
\vskip 5 pt

\indent Now, we let $\mathcal{K}$ be the set of all cliques $K_{k}$ of $G$. Moreover, we let  $\mathcal{K}_{1}$ be the set of all cliques $K_{k}$ of $G$ such that $V(K_{k}) \cap S \neq \emptyset$ and $\mathcal{K}_{2}$ be the set $\mathcal{K} \backslash \mathcal{K}_{1}$. Since $S$ is an $\iota_{k}$-set of $G$, every $K_{k} \in \mathcal{K}_{2}$ is adjacent to a vertex in $S$. We also let $\mathcal{K}_{3}$ be the subset of $\mathcal{K}_{2}$ such that all cliques $K_{k}$ of $\mathcal{K}_{3}$ are not adjacent to any vertex in $\{v_{1}, v_{2}, ..., v_{\ell}\}$. Since $S$ is a $\{K_k\}$-isolating set of $G$, every clique in $\mathcal{K}_{3}$ is adjacent to a vertex of $S \setminus \{v_1, v_2, \cdots , v_{\ell}\}$ and thus contains a vertex of $A$. Hence every clique in $\mathcal{K}_{3}$ is adjacent to $B$ where $B$ is an independent dominating set of $G[A]$. Therefore $\{v_1,v_2,\cdots ,v_{\ell}\}\cup B$ is an independent $\{K_k\}$-isolating set of $G$ and $i'_k(G)\le \ell +|B| \le \ell + |A|$. By Claim 3,
\begin{align}
\iota '_k(G) \le
\ell + \Sigma^{\ell}_{i = 1}deg_{S_{i - 1}}(v_{i})(\Delta - deg_{S_{i - 1}}(v_{i})).
\end{align}
\vskip 5 pt


\noindent Let deg$_{S_{i - 1}}(v_{i})=x_i$ and define two functions $f$ and $g : (\mathcal{R^{+}} \cup \{0\})^{\ell} \rightarrow \mathcal{R}$ by
\begin{align}
f(x_{1}, x_{2}, ..., x_{\ell}) = \ell + \Sigma^{\ell}_{i = 1}x_{i}(\Delta - x_{i})\notag
\end{align}
and
\begin{align}
g(x_{1}, x_{2}, ..., x_{\ell}) = x_{1} + x_{2} + \cdots + x_{\ell} - |S| + \ell .\notag
\end{align}
To find an upper bound on $\iota'_k(G)$, we look for the maximum of $f(x_{1}, x_{2}, ..., x_{\ell})$ under the condition, due to Claim 2,
$g(x_{1}, x_{2}, ..., x_{\ell})=0$.
Let
\begin{align}
F(x_{1}, x_{2}, ..., x_{\ell},\lambda)=f(x_{1}, x_{2}, ..., x_{\ell})-\lambda
g(x_{1}, x_{2}, ..., x_{\ell})
\end{align} with $\lambda \in \mathcal{R}.$
From the Lagrange's multipliers method, we get an extremum for $f$ by letting
\begin{align}
\frac{\partial F(x_{1}, ..., x_{\ell}, \lambda)}{\partial x_{i}} = \Delta - 2x_{i} - \lambda = 0\notag ~~~{\rm for ~ all}~~~ 1 \leq i \leq \ell
\end{align}
\noindent and
\begin{align}
\frac{\partial F(x_{1}, ..., x_{\ell}, \lambda)}{\partial \lambda} = x_{1} + x_{2} + \cdots + x_{\ell} - |S| + \ell = 0.\notag
\end{align}
For this extremum all the $x_i$'s are equal to $\frac{|S|}{\ell}-1$ and the extremum is equal to
\begin{align}M=\ell +(\frac{|S|}{\ell}-1)(\ell \Delta - |S|+\ell).
\end{align}
For the particular values $\{x_1,x_2,\cdots ,x_{\ell}\}=\{|S|-\ell , 0,\cdots, 0\}$ corresponding to the case $N_S(v_i)=N_S(v_1)$ for $2\le i\le \ell$, $f(x_1,\cdots ,x_{\ell})=\ell +(|S|-\ell)(\Delta -|S|+\ell)\le M.$ Therefore the extremum $M$ of $f$ is a maximum and
\begin{align}
\iota_{k}'(G] &\leq \ell + (\frac{\iota_{k}}{\ell} - 1)(\ell \Delta -\iota_{k}+\ell)\notag\\
           &=  -\frac{\iota_k^2(G)}{\ell} +\iota_k(G)(\Delta +2)-\ell \Delta.  \notag
\end{align}
\noindent This completes the proof.
\qed

\section{Proof of Theorem \ref{thm 2}}\label{sec 4}
We restate Theorem \ref{thm 2}.
\vskip 5 pt

\noindent \textbf{Theorem \ref{thm 2}} \emph{Let $G$ be a $K_{1, r}$-free graph where $r \geq 3$. Then $\iota'_{k} \leq (r - 2)(\iota_{k} - 1) + 1$.}\\
\proof Let $S$ be an $\iota_{k}$-set of $G$. Clearly, $|S| = \iota_{k}$. Let $I$ be a maximum independent set of $G[S]$. If $I = S$, then $S$ is an independent $\{K_{k}\}$-isolating set of $G$ implying that $\iota'_{k} \leq |S| = \iota_{k}$. This completes the proof because $r \geq 3$. Hence, we may assume that $S \setminus I \neq \emptyset.$ Let $A=N_G(S) \setminus N_G[I]$ and let $B$ be an independent dominating set  of $G[A]$. Consider a partition $B=\bigcup _{1\le i \le |S\setminus I|}B_i$ where $B_i\subseteq N_A(v_i)$ for each $v_i \in S\setminus I$ and $B_i \cap B_j = \emptyset$ for all $i$ and $j$. Note that some $B_i$ may be empty. For each nonempty $B_i$, $|B_i|\le r-2$ since $v_i$ has at least one neighbor in $I$ and $G$ is $K_{1,r}$-free. Therefore
\begin{align}
|B|\le \Sigma _{1\le i \le |S\setminus I|}|B_i|\le |S \setminus I|(r-2).
\end{align}

Since $S$ is a $\{K_k\}$-isolating set of $G$, every clique $K_k$ of $G$ non-adjacent to a vertex of $I$ has a vertex in $A$ and is thus adjacent to a vertex of $B$. Therefore $I\cup B$ is an independent $\{K_k\}$-isolating set of $G$ and
\begin{align}
\iota'_{k}(G) \le |I|+|B|\le |I|+(|S|-|I|)(r-2)
            =|S|(r-2)-|I|(r-3)\notag
\end{align}

\noindent which is maximized when $|I| = 1$. Hence,
\begin{align}
\iota'_{k}  &\leq (r - 2)\iota_{k} - (r - 3) = (r - 2)(\iota_{k} - 1) + 1\notag
\end{align}
which completes the proof.
\qed

\end{document}